\documentclass[11pt,amssymb]{amsart}
\usepackage{amssymb}
\usepackage{amsmath}
\usepackage[mathscr]{eucal}

\newcommand{\Z}{{\mathbb Z}}

\newcommand{\Br}{\mathrm{Br}}

\newtheorem{thm}{Theorem}[section]
\newtheorem{lemma}[thm]{Lemma}

\newtheorem{cor}[thm]{Corollary}

\newcommand{\gen}{\mathbf{gen}}

\newcommand{\uG}{\underline{G}}
\usepackage[all,cmtip]{xy}
.240pk scaled 1200 .240pk
\begin{document}

\title[Groups with good reduction and genus]{Algebraic groups with good reduction and the genus problem}

\author[I.~Rapinchuk]{Igor A. Rapinchuk}

\address{Department of Mathematics, Michigan State University, East Lansing, MI}

\email{rapinchu@msu.edu}

\begin{abstract}
We first provide an overview of several results dealing with the genus of a division algebra and highlight the role of ramification in its analysis. We then give a survey of recent developments on the genus problem for simple algebraic groups and its connections to the analysis of groups with good reduction.
\end{abstract}

\maketitle

\section{Introduction}\label{S:Introduction}

This paper is an expanded version of the author's talk at the Amitsur Centennial Symposium.
Our goal is to give an overview of some recent work on the surprising and exciting connections between the study of algebraic groups with good reduction and the {genus problem}, which is concerned with characterizing simple algebraic groups in terms of their maximal tori over the field of definition. These developments in fact comprise just one aspect of the emerging arithmetic theory of algebraic groups over higher-dimensional fields in which good reduction plays a central role --- we refer the reader to \cite{RR-Survey} for an in-depth account of good reduction for linear algebraic groups and its connections with various other directions, including local-global principles and weak-commensurability of Zariski-dense subgroups and applications to Riemann surfaces. The reader can also consult \cite{RR-Notices} for a more concise and less technical survey of these topics.

The initial focus of the genus problem was on the study of division algebras having the same maximal subfields. So, to establish the appropriate context, we will begin by discussing the genus of a division algebra and the role of ramification in its analysis. We will then give a precise statement of the genus problem for algebraic groups. In this setting, it appears that groups with good reduction are an adequate substitute for unramified division algebras. We will therefore present a brief discussion of good reduction and formulate our main Finiteness Conjecture for groups with good reduction over finitely generated fields. We will then indicate the connection between the Finiteness Conjecture and the genus problem. To conclude the paper, we will survey some of the available results.

\section{The genus problem for division algebras}\label{S-GenusDivAlg} Let $K$ be a field and $D_1$ and $D_2$ be two central division $K$-algebras of degree $n$. We say that $D_1$ and $D_2$ {\it have the same maximal subfields} if a degree $n$ field extension $P/K$ admits a $K$-embedding $P \hookrightarrow D_1$ if and only if it admits a $K$-embedding $P \hookrightarrow D_2$. Then one can ask the following natural question:

\vskip2mm

\noindent $(*)$ \parbox[t]{16cm}{\it Let $D_1$ and $D_2$ be central division algebras of the same degree. How are $D_1$ and $D_2$ related if they have the same maximal subfields?}

\vskip2mm

\noindent We note that this problem is similar in spirit to the following famous theorem of Amitsur \cite{Amitsur}:

\vskip2mm

\begin{thm}\label{T:Amitsur}{\rm (Amitsur)} Let $D_1$ and $D_2$ be finite-dimensional central division algebras over a field $K$ that have the same \emph{splitting fields}, i.e for a field extension $F$, the algebra $D_1 \otimes_K F$ is $F$-isomorphic to
a matrix algebra $\mathrm{M}_{n_1}(F)$ if and only if the algebra $D_2 \otimes_K F$ is isomorphic to a matrix algebra $\mathrm{M}_{n_2}(F)$. Then $n_1 = n_2$ and the classes $[D_1]$ and $[D_2]$ in the Brauer group $\mathrm{Br}(K)$ generate the same subgroup, $\langle [D_1] \rangle = \langle [D_2] \rangle$.

\end{thm}

The crucial point is that Amitsur's proof
relies in a very essential way on
{\it infinite} (non-algebraic) extensions of $K$ --- namely, so-called {\it generic splitting fields} (concrete examples of which are function fields of Severi-Brauer varieties). So, it is natural to ask if Amitsur's Theorem can be proved using only {\it finite} extensions of $K$. In other words, does the theorem's conclusion still hold if one only assumes that 
that $D_1$ and $D_2$ have the same {\it finite-dimensional} splitting fields or just the same maximal subfields? It turns out that the answer is (strongly) negative already over global fields. Indeed, this follows from the observation that, using 
the Albert-Brauer-Hasse-Noether theorem (see \cite[Ch. 18, \S18.4]{Pierce}), one can construct arbitrarily large collections of pairwise non-isomorphic cubic division algebras having the same maximal subfields over number fields (the same construction actually works for division algebras of any degree $d > 2$ --- cf. \cite[\S1]{CRR2}).
On the other hand, two quaternion division algebras over a number field that have the same quadratic subfields are necessarily isomorphic. Thus, even over number fields, question $(*)$ appears to be interesting.

Until about 10 years ago, no information at all was available on $(*)$ over any fields other than global. The following question along these lines was first asked in \cite[Remark 5.4]{PrRap-WC}:

\vskip2mm

\centerline{\parbox[t]{16cm}{\it Are quaternion division algebras over $\mathbb{Q}(x)$ determined uniquely up to isomorphism by their maximal subfields?}}

\vskip2mm

\noindent Shortly after it was formulated, this question was answered in the affirmative by D.~Saltman. In subsequent work, he and S.~Garibaldi \cite{GS} showed that the answer is still affirmative over the field of rational functions $k(x)$, where $k$ is any number field, and also in some other situations. This marked the starting point of the investigation of question $(*)$ over fields more general than global (we note
that a similar question, formulated in terms of finite-dimensional splitting fields, was considered in \cite{KrashMc}).

\vskip2mm

To quantify our discussion, it is convenient to introduce the notion of the {\it genus} of a division algebra (this terminology was suggested by L.H.~Rowen).


\vskip2mm

\noindent {\bf Definition 2.2.} {\it Let $D$ be a finite-dimensional central division algebra over a field $K$. Then the genus $\mathbf{gen}(D)$ of $D$ is defined to be the set of classes $[D'] \in \mathrm{Br}(K)$ represented by central division $K$-algebras $D'$ having the same maximal subfields as $D$.}

\vskip2mm

Broadly speaking, the general goal of the {\it genus problem} is to characterize the genus of a given division algebra --- note that this is essentially a reformulation of our original question $(*)$. Most recent work in this direction has dealt with the following two more precise questions:

\vskip2mm

\noindent $\bullet$ \parbox[t]{16cm}{{\it When does $\mathbf{gen}(D)$ reduce to a single element?} (Note that this is the case if and only if $D$ is determined uniquely up to isomorphism by its maximal subfields.)}

\vskip1mm

\noindent $\bullet$ \parbox[t]{16cm}{{\it When is $\mathbf{gen}(D)$ finite?}}

\vskip2mm

We should point out that over a number field $K$,
the Albert-Brauer-Hasse-Noether Theorem enables one to resolve both questions: namely, the genus of every quaternion division algebra is trivial (i.e., reduces to a single element), while the genus of any division algebra of higher degree is nontrivial but always finite (see \cite[Proposition 3.1]{CRR2} for the details).

Turning now to more general fields, we would first like to mention
the following theorem for rational function fields that was established in \cite{RR}.

\addtocounter{thm}{1}
\begin{thm}\label{T:A2}
{\rm (Stability Theorem)}
Assume that $\mathrm{char}\: k \neq 2$. If  $\vert \mathbf{gen}(\Delta) \vert = 1$  for any central division quaternion algebra $\Delta$ over $k$,
then  $\vert \mathbf{gen}(D) \vert = 1$  for any quaternion algebra $D$ over $k(x)$.
\end{thm}


\noindent (An analogous statement also holds for all division algebras having exponent two in the Brauer group --- cf. \cite{CRR01}.) Note that a consequence of Theorem \ref{T:A2} is that the genus of a quaternion algebra over the purely transcendental extension $k(x_1, \ldots , x_r)$ of a number field $k$ of any (finite) transcendence degree reduces to a single element. On the other hand, at this point, it is not known whether there exists a central quaternion division algebra $D$ over a finitely generated field $K$ of characteristic $\neq 2$ having nontrivial genus.

Next, observe that $\vert \mathbf{gen}(D) \vert > 1$ whenever $D$ does not have exponent two since in that case, the opposite algebra $D^{\small \mathrm{op}}$ is not isomorphic to $D$, but clearly has the same maximal subfields as $D$. We will therefore focus on finiteness properties of the genus. Let us first point out that over general fields, the genus $\mathbf{gen}(D)$ can be infinite. Indeed, adapting a construction that has been suggested by a number of people, including M.~Schacher, A.~Wadsworth, M.~Rost, S.~Garibaldi and D.~Saltman, J.~Meyer \cite{Meyer} produced examples of quaternion algebras over ``large" fields with infinite genus.
(By construction, these fields have infinite transcendence degree over the prime subfield.) Subsequently, S.~Tikhonov \cite{Tikh} extended this approach to construct examples of division algebras of any prime degree having infinite genus. On the other hand, over finitely generated fields, we have the following.



\begin{thm}\label{T:A3}
Let $K$ be a finitely generated field. Then for any finite-dimensional central division $K$-algebra $D$, the genus $\mathbf{gen}(D)$ is finite.
\end{thm}

There are two proofs of Theorem \ref{T:A3}, which can be found in \cite{CRR3} and \cite{CRR-Israel}. Both arguments depend on an analysis of ramification, but require somewhat different amounts of information about the unramified Brauer group.
Before indicating some of the main points, we first recall that if $F$ is a field equipped with a discrete valuation $v$, then a central simple $F$-algebra $A$ is said to be {\it unramified} at $v$ if there exists an Azumaya algebra $\mathcal{A}$ over the valuation ring $\mathcal{O}_v$ of the completion $K_v$ such that
$$
\mathcal{A} \otimes_{\mathcal{O}_v} K_v \simeq A \otimes_K K_v.
$$
Furthermore, if $F$ is equipped with a set $V$ of discrete valuations, the corresponding {\it unramified Brauer group} is defined as
$$
\Br(F)_V = \{ x \in \Br(F) \mid x \ \text{is unramified at all} \ v \in V \}.
$$
We also recall that any finitely generated field $K$ possesses natural sets of discrete valuations called {\it divisorial}. More precisely, let $\mathfrak{X}$ be a {\it model} of $K$, i.e. a normal separated irreducible scheme of finite type over $\Z$ (if $\mathrm{char}\: K = 0$) or over a finite field (if $\mathrm{char}\: K > 0$) such that $K$ is the function field of $\mathfrak{X}$. It is well-known that to every prime divisor $\mathfrak{Z}$ of $\mathfrak{X}$, there corresponds a discrete valuation $v_{\mathfrak{Z}}$ on $K$. Then
$$
V(\mathfrak{X}) = \{ v_{\mathfrak{Z}} \ \vert \ \mathfrak{Z} \ \ \text{prime divisor of} \ \ \mathfrak{X} \}
$$
is called the divisorial set of places of $K$ corresponding to the model $\mathfrak{X}$. Any set of places $V$ of $K$ of this form (for some model $\mathfrak{X}$) will be simply called {\it divisorial}.

The proof of Theorem \ref{T:A3} given in \cite{CRR3} requires the additional assumption that the degree $n$ of the division algebra $D$ is relatively prime to $\mathrm{char}\: K$. For the argument, we fix a divisorial set of places $V$ of $K$. Note that since $\mathrm{char}\: K$ is prime to $n$, we can assume without loss of generality that for each $v \in V$, the characteristic of the residue field $K^{(v)}$ is prime to $n$. Then one of the essential ingredients in the proof is the following observation:

\begin{lemma}\label{L:GenRam} If $D$ and $D'$ are central division $K$-algebras of degree $n$ having the same maximal subfields, then for any $v \in V$,
either both algebras are unramified at $v$  or both are ramified.
\end{lemma}

\noindent (See \cite[Lemma 2.5]{CRR01} for the details.) This ultimately enables one to reduce the proof of Theorem \ref{T:A3} to showing that the $n$-torsion subgroup $_n\Br(K)_V$ is finite and leads to the estimate
$$
\vert \gen(D) \vert \leq \vert {}_n\mathrm{Br}(K)_V \vert \cdot \varphi(n)^r,
$$
where $r$ is the number of $v \in V$ where $D$ ramifies (which is necessarily finite for a divisorial set).

Our second proof of Theorem \ref{T:A3}, given in \cite{CRR-Israel}, also uses the analysis of ramification, but does not impose any restrictions on the characteristic of $K$. The reason is that the argument does not require the finiteness of the full $n$-torsion subgroup of $\Br(K)$, but only the finiteness of certain of its subgroups. On the other hand,  since these subgroups depend on the division algebra at hand, we do not obtain a general estimate on the size of the genus as provided by our first proof.

Without going into the details, let us mention that an interesting generalization of the genus, termed the {\it upper genus}, was analyzed in \cite{KMRRS}, where the study of the genus problem for algebras with involution was also begun.

\section{The genus problem for algebraic groups} We now turn our attention to the genus problem for algebraic groups. Roughly speaking, to define the genus in this context, we replace maximal subfields with maximal tori in the definition of the genus of a division algebra.

More precisely, given two reductive algebraic groups $G_1$ and $G_2$ over a field $K$, we say that $G_1$ and $G_2$ {\it have the same isomorphism classes of maximal $K$-tori} if every maximal $K$-torus $T_1$ of $G_1$ is $K$-isomorphic to some maximal $K$-torus $T_2$ of $G_2$, and vice versa.
\vskip1mm

\noindent {\bf Definition 3.1.} {\it Let $G$ be an absolutely simple simply connected algebraic group over a field $K$. The {\it genus} $\mathbf{gen}_K(G)$ of $G$ is the set of $K$-isomorphism classes of (inner) $K$-forms $G'$ of $G$ that have the same isomorphism classes of maximal $K$-tori as $G$.}

\vskip1mm

\noindent (We recall that if $G$ is an algebraic group over $K$, then a $K$-group $G'$ is a called a {\it $K$-form of $G$} if $G$ and $G'$ become isomorphic over a separable closure $K^{\mathrm{sep}}$.)

\vskip1mm

In analogy with the case of division algebras, the following two questions have received the most attention so far.

\vskip1mm

\noindent $\bullet$ {\it When does $\mathbf{gen}_K(G)$ reduce to a single element?}

\vskip1mm

\noindent $\bullet$ {\it When is $\mathbf{gen}_K(G)$ finite?}

\vskip1mm

The basic case where $K$ is a number field was considered in \cite[Theorem 7.5]{PrRap-WC}, where the following result was established (although the term ``genus," which appeared later, was not used).

\addtocounter{thm}{1}

\begin{thm}\label{T:A1}
Let $G$ be an absolutely almost simple simply connected algebraic group over a number field $K$. Then

\vskip1.5mm

\noindent {\rm (1)} $\mathbf{gen}_K(G)$ is finite;

\vskip1.5mm

\noindent {\rm (2)} if $G$ is not of type $\textsf{A}_n$, $\textsf{D}_{2n+1}$ $(n > 1)$, or $\textsf{E}_6$, we have $\vert \mathbf{gen}_K(G) \vert = 1$.
\end{thm}

Given that this result resolves both questions over number fields, the next natural problem is to investigate the behavior of the genus over more general (finitely generated) fields. On the basis of Theorem \ref{T:A1} (as well a number of other results that we will discuss in \S\ref{S-Results}), in conjunction with the statements for division algebras mentioned in the previous section, we have been led to the following conjecture.

\vskip2mm

\noindent {\bf Conjecture 3.3.}

\vskip1mm

\noindent (1) \parbox[t]{16cm}{{\it Let $K = k(x)$ be the field of rational functions in one variable over a number field $k$. If $G$ is an absolutely almost simple simply connected algebraic $K$-group with center $Z(G)$ of order $\leq 2$, then the genus $\mathbf{gen}_K(G)$ reduces to a single element.}}

\vskip1mm

\noindent (2) \parbox[t]{16cm}{{\it Let $G$ be an absolutely almost simple simply connected algebraic group over a finitely generated field of ``good" characteristic. Then the genus
$\mathbf{gen}_K(G)$ is finite.}}

\vskip2mm

\noindent (Here, ${\rm char}~K = p$ is said to be ``good" if either $p = 0$ or $p >0$  and does not divide the order of the Weyl group of $G$.)

\vskip1mm

As we saw in the previous section, the study of the genus of a division algebra is based on a careful analysis of ramification. In the setting of algebraic groups, it appears that a suitable substitute for unramified division algebras are algebraic groups with good reduction. In the next section, we will describe the precise connection between the genus and groups with good reduction, which underlies most of the progress that has been achieved on Conjecture 3.3 so far.



\section{Groups with good reduction} We begin this section by first recalling some of the basic facts concerning groups with good reduction and then  formulating our main Finiteness Conjecture for forms with good reduction.

Let $K$ be a field equipped with a discrete valuation $v$, and suppose $G$ is a reductive affine algebraic group over $K$. We say that $G$ has {\it good reduction at $v$} if there exists a reductive
group scheme
\hspace{-2mm}\footnote{Let $R$ be a commutative ring and $S = \mathrm{Spec}\: R$. Recall that a {\it reductive $R$-group scheme} is a smooth affine group scheme $\mathscr{G} \to S$ such that the geometric fibers $\mathscr{G}_{\bar{s}}$ are connected reductive algebraic groups (see \cite[Exp. XIX, Definition 2.7]{DemGr} or \cite[Definition 3.1.1]{Con}).}
$\mathscr{G}$ over the valuation ring $\mathcal{O}_v$ of the completion $K_v$ whose generic fiber $\mathscr{G} \otimes_{\mathcal{O}_v} K_v$ is isomorphic to $G \otimes_K K_v$. The special fiber (or {\it reduction}) $\mathscr{G} \otimes_{\mathcal{O}_v} K^{(v)}$,  where $K^{(v)}$ is the residue field of $K_v$, is then denoted by $\underline{G}^{(v)}$ --- it is a connected algebraic group over $K^{(v)}$ of the same type as $G$.
Furthermore, given a set $V$ of discrete valuations of $K$, we will say that $G$ has {\it good reduction with respect to $V$} if $G$ has good reduction at all $v \in V$.

Informally, having good reduction means that the group $G \times_K K_v$ has a ``nice" $\mathcal{O}_v$-structure whose reduction modulo $\mathfrak{p}_v$ yields a connected reductive group. As the following examples demonstrate, in various situations of interest, this condition can
be characterized in very concrete terms.

\vskip3mm

\noindent {\bf Example 4.1.}


\vskip1mm

\noindent (a) \parbox[t]{16cm}{(cf. \cite[Example 2.2]{CRR-AlgGpGenus}) If $G = \mathrm{SL}_{1,A}$ is the algebraic group associated with the elements of reduced norm 1 in a central simple algebra $A$ over $K$, then $G$ has good reduction at $v$ if and only if $A$ is unramified at $v$.}

\vskip1mm

\noindent (b) \parbox[t]{16cm}{(cf. \cite[Example 2.3]{CRR-AlgGpGenus}) Suppose $q$ is a non-degenerate quadratic form in $n$ variables over $K$ and the residue field $K^{(v)}$ has characteristic $\neq 2.$ Then the spinor group $G = {\rm Spin}_n (q)$ has good reduction at $v$ if and only if $q$ is equivalent over $K_v$ to a quadratic form
$$
\lambda(a_1 x_1^2 + \cdots + a_n x_n^2), \ \ \ \text{with} \ \lambda \in K_v^{\times}, a_i \in \mathcal{O}_v^{\times}.
$$}

\vskip2mm

Historically, the study of algebraic groups with good reduction can be traced back to the work of Harder \cite{Harder}, Colliot-Th\'el\`ene and Sansuc \cite{CTS}, and Gross \cite{Gross}. In \cite{Harder}, the focus was mainly on algebraic groups over a number field $K$ having good reduction with respect to sets $V$ consisting of almost all nonarchimedean places of $K$. Groups with good reduction have also been analyzed extensively over $K = k(x)$ (the field of rational functions in one variable over a field $k$). In \cite{RagRam}, Raghunathan and Ramanathan considered the case where $V$ consists of the discrete valuations $v_{p(x)}$ corresponding to {\it all} irreducible polynomials $p(x) \in k[x]$. Later, groups having good reduction at all $v \in V \setminus \{ v_x \}$ were analyzed by Chernousov, Gille and Pianzola \cite{CGP}; these results then played a crucial role in their proof of the conjugacy of analogues of Cartan subalgebras in certain infinite-dimensional Lie algebras.

A key point is that in all of these cases, $K$ is the fraction field of a Dedekind ring $R$, and $V$ consists of discrete valuations associated with the nonzero prime ideals of $R$, making the situation ``1-dimensional." By contrast, our recent and ongoing work addresses the analysis of groups with good reduction in the higher-dimensional setting of arbitrary finitely generated fields. More precisely, the following Finiteness Conjecture is one of the central elements of our current investigations.

\vskip2mm

\noindent \noindent {\bf Conjecture 4.2.} {\it Let $G$ be a reductive algebraic group over a finitely generated field $K$ and $V$ be a divisorial set of discrete valuations of $K$. Then the set of $K$-isomorphism classes of $K$-forms $G'$ of $G$ that have good reduction at all $v \in V$ is finite (at least when the characteristic of $K$ is ``good.").}

\vskip2mm

\noindent (When $G$ is an absolutely almost simple algebraic group, ``good" characteristic is used here in the same sense as in Conjecture 3.3. For algebraic tori, by good characteristic, we simply mean ${\rm char}~K = 0.$) 

\vskip2mm

We refer the reader to \cite{RR-Survey} for an extensive discussion of groups with good reduction and, in particular, the key role of the Finiteness Conjecture in the current development of the arithmetic theory of algebraic groups over higher-dimensional fields. For our present purposes, we will only mention some connections to the genus problem.

\addtocounter{thm}{2}

\begin{thm}\label{T:GoodReduction1} {\rm (\cite[Theorem 5]{CRR-Genus}, \cite[Theorem 1.1]{CRR-AlgGpGenus})}
Let $G$ be an absolutely almost simple linear algebraic group over a field $K$ and let $v$ be a discrete valuation of $K$. Assume that the residue field $K^{(v)}$ is finitely generated and that $\mathrm{char}\: K^{(v)} \neq 2$ if $G$ is of type $\textsf{B}_{\ell}$ $(\ell \geq 2)$.
If $G$ has good reduction at $v$, then any $G' \in \gen_K(G)$ also has good reduction at $v$. Moreover, the
reduction ${\uG'}^{(v)}$ lies in the genus $\gen_{K^{(v)}}(\uG^{(v)})$ of the reduction $\uG^{(v)}$.
\end{thm}

As we already mentioned above, one should view groups with good reduction as an analogue of unramified division algebras (this point of view is justified, for instance, by Example 4.1). From this perspective, Theorem \ref{T:GoodReduction1} can be thought of as a partial analogue of Lemma \ref{L:GenRam}. Although we refer the reader to \cite{CRR-AlgGpGenus} for the details, we would like to point out that the proof of Theorem \ref{T:GoodReduction1} is based on an entirely new approach to good reduction of simple algebraic groups that shows
that the existence of good reduction can be characterized in terms of the presence of maximal tori with certain specific properties.

In the case of finitely generated fields, we have the following consequence.

\begin{cor}\label{C:genus-GR}
Let $G$ be an absolutely almost simple algebraic group over an infinite finitely generated field $K$, and let $V$ be a divisorial set of places of $K$. Assume that $\mathrm{char}\: K \neq 2$ if $G$ is of type $\textsf{B}_{\ell}$ $(\ell \geq 2)$. Then there exists a finite subset $S \subset V$ such that every $G' \in \gen_K(G)$ has good reduction at all $v \in V \setminus S$.
\end{cor}

In particular, it follows that the truth of Conjecture 4.2 for all divisorial sets $V$ would yield the finiteness of $\gen_K(K)$ for any absolutely almost simple algebraic $K$-group $G$. In other words, the Finiteness Conjecture provides a uniform approach for resolving one important aspect of the genus problem for algebraic groups. Continuing the parallel with division algebras, we thus see that, in this context, the Finiteness Conjecture plays a role analogous to that of our finiteness results for the unramified Brauer group in the study of the genus of a division algebra.

\section{Overview of results}\label{S-Results} In this section, we will give a brief overview of some available results on Conjecture 4.2 and the genus problem. The reader can consult \cite{RR-Survey} for a more detailed account.

To begin with, Conjecture 4.2 has been settled completely for algebraic tori.



\begin{thm}\label{T-ToriGR} {\rm (cf. \cite[Theorem 1.1]{RR-Tori}, \cite[Corollary 5.2]{RR-Tori1})} Let $K$ be a finitely generated field and $V$ be a divisorial set of places of $K$. Then for any $d \geq 1$, the set of $K$-isomorphism classes of $d$-dimensional $K$-tori that have good reduction at all $v \in V$ and for which the degree $[K_T : K]$ of the minimal splitting field is prime to the characteristic exponent of $K$, is finite.
\end{thm}

Turning now to semisimple groups, we first note that in the classical setting when $K$ is a number field and $V$ is any set consisting of almost all nonarchmidean places of $K$, Conjecture 4.2 can be reduced to the consideration of absolutely almost simple groups, in which case the assertion
follows from well-known results on the Galois cohomology of algebraic groups and the description of groups with good reduction over $p$-adic fields (see \cite[Proposition 5.2]{RR-Survey}). The general case, however, presents a number of new challenges.

Here are some representative results. First, we have the following statement for inner forms of type $\mathsf{A}_n$ over (essentially arbitrary) finitely generated fields, which is derived from our finiteness results for the unramified Brauer group (discussed in \S\ref{S-GenusDivAlg}).

\begin{thm}\label{T-GRAn} {\rm (\cite{CRR3})}
Let $K$ be a finitely generated field, $V$ a divisorial set of discrete valuations of $K$, and $n \geq 2$ an integer prime to $\mathrm{char}\: K$. Then the number of $K$-isomorphism classes of groups of the form $\mathrm{SL}_{1,A}$, where $A$ is a central simple $K$-algebra of degree $n$, that have good reduction at all $v \in V$, is finite.
\end{thm}

Using this statement, together with Theorems \ref{T:A2} and \ref{T:GoodReduction1} and some additional considerations involving generic tori, we obtain the next result concerning the genus.

\begin{thm}\label{T:A5} {\rm (cf. \cite[Theorem 5.3]{CRR01} and \cite[Theorem 6.3]{CRR2})}
\vskip1mm

\noindent {\rm (1)} \parbox[t]{16cm}{Let $D$ be a central division algebra of exponent two over the field of rational functions $K = k(x_1, \ldots , x_r)$, where $k$ is
either a number field or a finite field of characteristic $\neq$ 2. Then for $G = \mathrm{SL}_{m , D}$ $(m \geq 1)$, the genus $\mathbf{gen}_K(G)$
reduces to a single element.}

\vskip1mm

\noindent {\rm (2)} \parbox[t]{16cm}{Let $G = \mathrm{SL}_{m , D}$, where $D$ is a central division algebra over a finitely generated field $K$ of degree prime to $\mathrm{char}\: K$. Then $\mathbf{gen}_K(G)$ is finite.}
\end{thm}

Next, following Kato \cite{Kato}, we say that a field $K$ is a {\it 2-dimensional global field} if it is either the function field $k(C)$
of a smooth geometrically integral curve $C$ over a number field $k$ or the funcion field $\mathbb{F}_q(S)$ of a smooth geometrically integral surface $S$ over a finite field $\mathbb{F}_q.$

\begin{thm}\label{T-GR} {\rm (cf. \cite[Theorem 1.1]{CRR4})}
Let $K$ be a two-dimensional global field of characteristic $\neq 2$ and let $V$ be a divisorial set of discrete valuations of $K$. Fix an integer $n \geq 5$. Then the set of $K$-isomorphism classes of spinor groups $G = \mathrm{Spin}_n(q)$ of nondegenerate quadratic forms in $n$ variables over $K$ that have good reduction at all $v \in V$ is finite.
\end{thm}

Using Voevodksy's resolution of Milnor's conjecture on quadratic forms, we reduce the proof of this result to the analysis of finiteness properties of unramified cohomology with $\mu_2$-coefficients, which we verify in the case of 2-dimensional global fields. Similar statements are also available for some special unitary groups of types $\textsf{A}_n$ and $\textsf{C}_n$ and for groups of type $\textsf{G}_2$. We note that if $n$ is odd, then all $K$-forms of $G = {\rm Spin}_n(q)$ are of the form $G' = {\rm Spin}_n(q')$, with $q'$ a non-degenerate quadratic form over $K$ in $n$ variables, Theorem \ref{T-GR} yields Conjecture 4.2 in this case.

Turning now to the genus, we have the following statements for spinor groups, groups of type $\mathsf{G}_2$, and groups of type $\mathsf{F}_4$ that split over a quadratic extension.

\begin{thm}\label{T:A6} {\rm (cf. \cite[Theorem 1.2]{CRR4}, \cite[Theorem 5.5]{RR-Tori})}
Suppose $K$ is either a 2-dimensional global field of characteristic $\neq 2$ or the field of rational functions $k(x , y)$ in two variables over a number field $k$. Let $G = \mathrm{Spin}_n(q)$ be the spinor group of a nondegenerate quadratic form $q$ over $K$ of \emph{odd} dimension $n \geq 5$. Then $\mathbf{gen}_K(G)$ is finite.
\end{thm}

\begin{thm}\label{T:A7} {\rm (cf. \cite[Theorems 9.1 and 9.3]{CRR4} and \cite[Proposition 5.3]{RR-Tori})}
Let $G$ be a simple algebraic $K$-group of type $\textsf{G}_2$.

\vskip2mm

\noindent {\rm (1)} If $K$ is the field of rational functions $k(x)$, where $k$ is a number field, then $\vert \mathbf{gen}_K(G) \vert = 1$.

\vskip2mm

\noindent {\rm (2)} \parbox[t]{16cm}{If $k$ is a number field and $K$ is one of the following:

\vskip2mm

\hskip5mm $\bullet$ $K = k(x_1, \ldots, x_r)$ is the field of rational functions in any (finite) number of variables;

\vskip1mm

\hskip5mm $\bullet$ $K = k(C)$ is the function field of a smooth geometrically integral curve $C$ over $k$;

\vskip1mm

\hskip5mm $\bullet$ \parbox[t]{15cm}{$K = k(X)$ is the function field of a Severi-Brauer variety $X$ over $k$ associated with a central division algebra $D$ over $k$ of degree $\ell$, where $\ell$ is either odd or $\ell = 2$,}

\vskip2mm

\noindent then $\mathbf{gen}_K(G)$ is finite.}
\end{thm}

\begin{thm}\label{T:F4GR} {\rm (cf. \cite[Theorems 1.10 and 1.11]{CRR-AlgGpGenus})} \ \

\vskip1mm

\noindent {\rm (1)} \parbox[t]{16cm}{Let $k$ be a number field, and set $K = k(x)$. Then for any simple algebraic $k$-group $G$ of type $\textsf{F}_4$ that splits over a quadratic extension of $K$, the genus $\gen_K(G)$ is trivial.}

\vskip1mm

\noindent {\rm (2)} \parbox[t]{16cm}{Let $K$ be either a 2-dimensional global field of characteristic $\neq 2, 3$ or the field of rational functions $k(x,y)$ in two variables over a number field $k$. Then for any simple $K$-group $G$ of type $\textsf{F}_4$ that splits over a quadratic extension of $K$, the genus $\gen_K(G)$ is finite.}

\end{thm}

To conclude this section, we would like to mention a couple of results dealing with a newly-discovered phenomenon that we refer to as ``killing the genus by a purely transcendental extension" --- the reason for this choice of terminology is that, in the two cases considered below, the genus essentially becomes as small as possible after passing to a suitable purely transcendental extension.

\begin{thm}\label{T:Kill1}{\rm (cf. \cite[Theorem 1.5]{CRR-AlgGpGenus})}
Let $A$ be a central simple algebra of degree $n$ over a finitely generated field $k$, and let $G = \mathrm{SL}_{1 , A}$. Assume that $\mathrm{char}\: k$ is prime to $n$, and let $K = k(x_1, \ldots , x_{n-1})$ be the field of rational functions in $(n-1)$ variables. Then $\gen_K(G \times_k K)$ consists of (the isomorphism classes of) groups of the form $H \times_k K$, where $H = \mathrm{SL}_{1 , B}$ and $B$ is a central simple algebra of degree $n$ such that its class $[B]$ in the Brauer group $\Br(k)$ generates the same subgroup as the class $[A]$.
\end{thm}

The proof uses Amitsur's theorem on generic splitting fields \cite{Amitsur}, and a result of D.~Saltman \cite{Salt1}, \cite{Salt2} on
function fields of Severi-Brauer varieties.

\begin{thm}\label{T:Kill2}{\rm (cf. \cite[Theorem 1.6]{CRR-AlgGpGenus})}
Let $G$ be a group of type $\textsf{G}_2$ over a finitely generated field $k$ of characteristic $\neq 2, 3$, and let $K = k(x_1, \ldots , x_6)$ be the field of rational functions in 6 variables. Then $\gen_K(G \times_k K)$ reduces to a single element.
\end{thm}

The proof relies on properties of Pfister forms.

\vskip5mm

\noindent {\small {\bf Acknowledgements.} I would like to thank the conference organizers for the opportunity to give a talk at the Amitsur Centennial Symposium. I would also like to thank the anonymous referee for a careful reading of the paper and suggestions that helped to improve the exposition. I was partially supported by a Collaboration Grant for Mathematicians from the Simons Foundation and by NSF grant DMS-2154408.}

\vskip5mm

\bibliographystyle{amsplain}

\begin{thebibliography}{100}

\bibitem{Amitsur} S.A.~Amitsur, {\it Generic splitting fields of central simple algebras}, Ann. of math. (2) {\bf 62}(1955), 8-43.

\bibitem{CGP} V.I.~Chernousov, P.~Gille, A.~Pianzola, {\it Torsors over the punctured affine line}, Amer. J. Math. {\bf 134} (2012), no. 6, 1541-1583.

\bibitem{CRR01} V.I.~Chernousov, A.S.~Rapinchuk, and I.A.~Rapinchuk, {\it The genus of a division algebra
and the unramified Brauer group},  Bull. Math. Sci. {\bf 3} (2013), 211-240.

\bibitem{CRR2}  V.I.~Chernousov, A.S.~Rapinchuk, I.A.~Rapinchuk, {\it Division algebras with the same maximal subfields}, Russian Math. Surveys {\bf 70}(2015), no. 1, 91-122.

\bibitem{CRR3}  V.I.~Chernousov, A.S.~Rapinchuk, I.A.~Rapinchuk, {\it On the size of the genus of a division algebra,} Proc. Steklov Inst. of Math. {\bf 292}(2016), no. 1, 63-93.

\bibitem{CRR-Genus} V.I.~Chernousov, A.S.~Rapinchuk, and I.A.~Rapinchuk, {\it On some finiteness properties of algebraic groups over finitely generated fields}, C. R. Acad. Sci. Paris, Ser. I {\bf 354} (2016), 869-873.


\bibitem{CRR4} V.I.~Chernousov, A.S.~Rapinchuk, and I.A.~Rapinchuk, {\it Spinor groups with good reduction}, Compos. Math. {\bf 155} (2019), no. 3, 484-527.

\bibitem{CRR-Israel} V.I.~Chernousov, A.S.~Rapinchuk, I.A.~Rapinchuk, {\it The finiteness of the genus of a finite-dimensional division algebra, and generalizations}, Israel J. Math. {\bf 236} (2020), no. 2, 747-799.

\bibitem{CRR-AlgGpGenus} V.I.~Chernousov, A.S.~Rapinchuk, I.A.~Rapinchuk, {\it Simple algebraic groups with the same maximal tori, weakly commensurable Zariski-dense subgroups, and good reduction},  arXiv:2112.04315

\bibitem{CTS} J.-L.~Colliot-Th\'el\`ene, J.-J.~Sansuc, {\it Fibr\'es quadratiques et composantes connexes r\'eelles}, Math. Ann. {\bf 244} (1979), no. 2, 105-134.

\bibitem{Con} B.~Conrad, {\it Reductive group schemes}, Autour des sch\'emas en groupes, \'Ecole d'\'et\'e ``Sch\'emas en groupes," vol. I (Luminy 2011), Soc. Math. France, Paris 2014.

\bibitem{DemGr} M.~Demazure, A.~Grothendieck, {\it Sch\'emas en groupes},
S\'eminaire
de g\'eom\'etrie alg\'ebrique du Bois Marie 1962/64 (SGA 3). Lect. Notes Math. {\bf 151}, {\bf 152}, {\bf 153}, Springer, 1970.

\bibitem{GS} S.~Garibaldi, D.~Saltman, {\it Quaternion Algebras
with the Same Subfields}, Quadratic forms, linear algebraic groups,
and cohomology,  225-238, Dev. math. 18, Springer, New York, 2010.

\bibitem{Gross} B.H.~Gross, {\it Groups over $\Z$}, Invent. math. {\bf 124}(1996), no. 1-3, 263-279.

\bibitem{Harder} G.~Harder, {\it Halbeinfache Gruppenschemata \"uber Dedekindringen}, Invent. math. {\bf 4}(1967), 165-191.


\bibitem{Kato} K.~Kato, {\it A Hasse principle for two dimensional global fields, with an appendix by J.-L.~Colliot-Th\'el\`ene}, J. reine angew. Math.
{\bf 366}(1986), 142-183.

\bibitem{KrashMc} D.~Krashen, K.~McKinnie, {\it Distinguishing algebras by their finite splitting fields}, Manuscr. math. {\bf 134} (2011), no. 1-2, 171-182.

\bibitem{KMRRS} D.~Krashen, E.~Matzri, A.S.~Rapinchuk, L.H.~Rowen, D.~Saltman, {\it Division algebras with common subfields}, Manuscr. math., https://doi.org/10.1007/s00229-021-01315-5

\bibitem{Meyer} J.S.~Meyer, {\it A division algebra with infinite genus}, Bull. London Math. Soc. {\bf 46}(2014), 463-468.


\bibitem{Pierce} R.S.~Pierce, {\it Associative Algebras}, GTM 88, Springer, 1982.

\bibitem{PrRap-WC} G.~Prasad, A.S.~Rapinchuk, {\it Weakly commensurable arithmetic
groups and isospectral locally symmetric spaces},
Publ. math. IHES {\bf 109}(2009), 113-184.

\bibitem{RagRam} M.S.~Raghunathan, A.~Ramanathan, {\it Principal bundles on the affine line}, Proc. Indian Acad. Sci. (Math. Sci.) {\bf 93}(1984), nos. 2-3, 137-145.

\bibitem{RR} A.S.~Rapinchuk,  I.A.~Rapinchuk, {\it On division algebras having the
same maximal subfields}, Manuscr. math. {\bf
132} (2010), 273-293.

\bibitem{RR-Survey} A.S.~Rapinchuk, I.A.~Rapinchuk, {\it Linear algebraic groups with good reduction}, Res. Math. Sci. {\bf 7} (2020), no. 3, 28.

\bibitem{RR-Notices}  A.S.~Rapinchuk, I.A.~Rapinchuk, {\it Recent developments in the theory of linear algebraic groups: Good reduction and finiteness properties}, Notices Amer. Math. Soc. {\bf 68} (2021), no. 6, 899-910.

\bibitem{RR-Tori} A.S.~Rapinchuk, I.A.~Rapinchuk, {\it Some finiteness results for algebraic groups and unramified cohomology over higher-dimensional fields}, J. Number Theory {\bf 233} (2022), 228-260.

\bibitem{RR-Tori1} A.S.~Rapinchuk, I.A.~Rapinchuk, {\it Properness of the global-to-local map for algebraic groups with toric connected component and other finiteness properties}, to appear in Mathematical Research Letters.

\bibitem{Salt1} D.~Saltman, {\it Norm polynomials and algebras}, J. Algebra {\bf 62}(1980), 333-345.


\bibitem{Salt2} D.~Saltman, {\it Lectures on Division Algebras}, CBMS Regional Conference Series in Mathematics {\bf 94}, AMS 1999.


\bibitem{Tikh} S.V.~Tikhonov, {\it Division algebras of prime degree with infinite genus}, Proc. Steklov Inst. {\bf 292} (2016), 256-259.

\end{thebibliography}

\end{document}